# A Problem of W. R. Scott: Classify the Subgroup of Elements with Many Roots

*Vance Faber*

**Abstract.** Let *G* be an infinite group and let *h* and *g* be elements. We say that *h* is a root of *g* if some integer power of *h* is equal to *g*. We define *K*(*G*) to be the subgroup of all elements of *G* for which the number of elements which are not roots is of smaller cardinality than the cardinality of the group. That is, each element in *K* has almost every element in *G* as a root. This paper discusses the problem: When can *K*(*G*) be non-trivial?

**§1. Introduction**. Nearly 60 years ago in [8], W. R. Scott defined two subgroups of a given infinite group *G*. The first, $K = K(G)$, is defined as follows:

For $g \in G$, let $\eta(g, G)$ be the set of $h \in G$ such that $h^n = g$ has no solution for $n$. (In other words, $\eta(g) = \{h : g \notin <h>\}$.) Then

$$K(G) = \{k \in G : |\eta(k)| < |G|\}.$$

The second group, $D = D(G)$, is defined by

$$D = \cap \{H \leq G : |H| = |G|\}.$$

Scott proved several interesting theorems about *D* and *K*; the principal ones follow.

**THEOREM A.** *If G is abelian, $D = K = E$ unless $G \cong Z_{p^\infty} \times F$, F finite, in which case, $D = K \cong Z_{p^\infty}$.*

**THEOREM B.** *For any G, $K \leq D$. If G is not periodic, $K = E$.*

**THEOREM C.** *For any G,*

*(i) $K \leq Z(G)$.*

*(ii) K is either cyclic of order $p^n$ or a $p^\infty$-group for some prime p.*

*(iii) K is a $p^\infty$-group if and only if there exists a central $p^\infty$-subgroup C such that $G/C$ is finite. If such a C exists, then $C = K = D$.*

We attack two main questions in this paper.

Question 1. Are there locally finite groups with $D > K$?

Question 2. What are the groups with non-trivial $K$?

We answer the first question affirmatively in Example 4.10.

Question 2 seems much more difficult to answer. For countable binary finite groups we have Corollary 4.6.

**COROLLARY 4.6.** *If G is countable and either binary finite or a 2-group such that $K(G) \neq E$, then either (i) $G = (P \times F)_C$, with F a finite group, P a $p^\infty$-group, and $C \leq Z(F)$ a cyclic group of order a power of p, or (ii) $G = \langle x, G_1 \rangle$, where $G_1 = (P \times F)_C$ with F a finite group, P a $2^\infty$-group, $C \leq Z(F)$ a cyclic group of order a power of 2; $x^{-1}zx = z^{-1}$ for every element z of P, $x^2 \in G_1$ and there exists an m such that $x^{2m} = a$, the unique element of P of order 2. In the first case, $K(G) = P$, and in the second case, $K(G) = \langle a \rangle$.*

For uncountable groups, we have Corollary 4.9.

**COROLLARY 4.9.** *If G is an uncountable class 2 nilpotent group, then $K(G) = E$.*

Clearly in Jonsson groups [**9**] $G$, $D(G) = G$. Thus, in discussions of $D$ we must put some restrictions on the class of groups we study, and locally finite seems a reasonable restriction since there are no uncountable locally finite Jonsson groups of regular cardinality (see [**2, p.73**]). However, the groups in Example 4.10 are nilpotent of class 2 and are even FC.

If G has an equipotent abelian subgroup, then by Theorem A, we know just about all there is to know about D and K, so we must expect to consider groups without equipotent abelian subgroups.

This paper is arranged as follows. Section 2 provides the notation. Section 3 gives a number of lemmas which are needed later. Section 4 provides the answer to Question 1 and concentrates on reducing Question 2 to seemingly more manageable questions.

**§2. Notation.** Let *S* and *T* be sets. The *cardinality* (*power*) of S is denoted by $|S|$. If *S* and *T* have the same cardinality, we say that they are equipotent. If *m* is an infinite cardinal, $m^+$ is the first cardinal greater than *m*; $2^m$ is the cardinality of the set of all subsets of a set of power *m*; $2^{<m} = \sum_{n<m} 2^n$; *m* is *regular* if it is not the sum of a smaller number of smaller cardinals.

Let *G* be a group, *H* be a subgroup, and *S* be a subset of *G*. The *centralizer* of *S* in *G* is denoted by $C_G(S)$. We denote the *center* of *G* by $Z(G) = C_G(G)$. If $H = \{1\}$, the trivial subgroup, we often denote *H* by *E*. If $x, y \in G$, the *commutator* of *x* and *y* is $[x, y] = x^{-1}y^{-1}xy$. The subgroup generated by $S$ is denoted by $<S>$. The *derived* (*commutator*) group of $G$ is $G' = <\{[x, y] : x, y \in G\}>$. The *index* of $H$ in $G$ is denoted by $[G : H]$. The *conjugate class* of $S$ in $G$ is $Cl(S) = Cl_G(S) = \{g^{-1}Sg : g \in G\}$. An element $g$ of $G$ is a $p'$-element if its order, $|g|$, is relatively prime to the prime $p$. The *exponent* of $G$ is the smallest integer $n$ (if one exists) such that $g^n = 1$ for all $g \in G$; it is denoted $\exp(G)$. $G$ is *nilpotent* of class $n$ if $n$ is the length of the upper central series. A *section* of $G$ is a factor group of a subgroup. $G$ is *FC* if every element $y$ has finitely many distinct conjugates $g^{-1}yg$ for $g \in G$. If $G$ is a $p$-group, $\Omega_1(G)$ is the group generated by the set of all elements of order $p$. A *Jonsson group* is an infinite group which has no proper equipotent subgroups. $G$ is *locally* (*binary*) *finite* if every finite (two-element) subset generates a finite subgroup. The cartesian product of $H$ and $C \leq G$ with amalgamated subgroup $A$ is denoted by $(H \times C)_A$. Additional terminology can be found in [**6**] and [**7**].

We also have occasion to use the following construction. Let $V$ and $W$ be vector spaces over a field $K$ and $\rho : V \times V \to W$ be a nonzero alternating bilinear function. If $\gamma : V \times V \to W$ is any bilinear function such that $\rho(x, y) = \gamma(x, y) - \gamma(y, x)$, then $V \times W$ can be given the structure of a nilpotent group of class 2, denoted by $G = V\gamma W$, by defining $(x, a)(y, b) = (x + y, a + b + \gamma(x, y))$. Note that $[(x, a), (y, b)] = (0, \rho(x, y))$.

**§3. Preliminary Lemmas**. Throughout this section $G$ is a group and $H$ is a subgroup. If $H \triangleleft G$ and $S \leq G$, we denote the set $SH/H$ by $\bar{S}$.

LEMMA 3.1. [**8, p.189**]

(i) $\eta(1) = \varnothing$.

(ii) $\eta(x^{-1}) = \eta(x)$.

(iii) $\eta(x_1 x_2) \leq \eta(x_1) \cup \eta(x_2)$.

(iv) $\eta(\sigma(x), \sigma(G)) \leq \sigma(\eta(x, G))$ *for any homomorphism* $\sigma$ *of G.*

(v) $\eta(h, H) = H \cap \eta(h, G)$ *if* $h \in H \leq G$.

(vi) $\eta(x) \geq G - C_G(x)$.

LEMMA 3.2.

(i) $x_1 \notin \eta(x_2)$ *and* $x_2 \notin \eta(x_3)$ *imply that* $x_1 \notin \eta(x_3)$.

(ii) $g \notin \eta(y)$ *and* $g \in \eta(a)$ *imply that* $y \in \eta(a)$.

PROOF.

(i) If $x_1 \notin \eta(x_2)$, there exists an $n$ such that $x_2 = x_1^n$. If $x_2 \notin \eta(x_3)$, there exists an $m$ such that $x_3 = x_2^m$. Thus $x_3 = x_2^{nm}$, so $x_1 \notin \eta(x_3)$.

Of course, (ii) is a contrapositive of (i).

LEMMA 3.3. Let $K = <a>$, $\eta(a) = \{g : a \notin <g>\}$. If $\alpha \in Aut(G)$, then $\alpha(\eta(a)) = \eta(a)$.

PROOF. We know by [**8, p.193**] that $K$ is a characteristic subgroup. Suppose $g \in \eta(a)$. If $[\alpha(g)]^n = a$, then $\alpha(g^n) = a$, so $g^n = \alpha^{-1}(a) \in K$. Since $\alpha(K) = K$, $\alpha^{-1}(a)$ is a generator of $K$. Thus there exists $m$ such that $(g^n)^m = (\alpha^{-1}(a))^m = a$, so $g \notin \eta(a)$. This contradiction shows that $\alpha(g) \in \eta(a)$. Since $\alpha^{-1}$ is also an automorphism, $\alpha^{-1}(\eta(a)) \leq \eta(a)$.

LEMMA 3.4. *If* $H \triangleleft G$, $|G/H| = |G|$ *and* $a \in K$, *then* $\bar{a} \in K(G/H)$.

PROOF. This follows from Lemma 3.1.(iv).

Next we generalize [**8; Corollary 3**].

LEMMA 3.5. *If* $H \leq K(G)$ *with $H$ finite, then* $\bar{K} = K(\bar{G})$.

PROOF. $\bar{K} \leq K(\bar{G})$ by Lemma 3.4. Suppose $\bar{x} \in K(\bar{G})$. Then $|\eta(\bar{x})| < |\bar{G}|$. Now

$$\eta(\bar{x}) = \{\bar{g} : \bar{x} \notin <\bar{g}>\} = \{\bar{g} : x \notin <g>H\}.$$

Let $\eta = \{g : \bar{g} \in \eta(\bar{x})\} = \{g : x \notin <g>H\}$. Then $|\eta| = |H||\eta(\bar{x})| < |G|$. If $g \notin \eta \cup (\bigcup_{h \in H} \eta(h))$, then $g^n = hx$ for some $n$ and some $h \in H$ and also $g^m = h^{-1}$ for some $m$. Hence $g^{m+n} = x$, so $g \notin \eta(x)$. Thus

$$|\eta(x)| \leq |\eta| + \sum_{h \in H} |\eta(h)| < |G|,$$

so $x \in K$.

LEMMA 3.6. *Suppose* $G = C(\eta(a))$, $1 \neq a \in K$, $|a| = p^n$. *Then the set $H$ of $p'$-elements is an abelian subgroup and* $H \leq \eta(a)$.

PROOF. Clearly if $x$ is a $p'$-element, $a \notin <x>$. Since $\eta(a)$ is abelian, the $p'$-elements form an abelian group.

LEMMA 3.7. *Suppose* $G = C(\eta(a))$, $1 \neq a \in K$, $|a| = p^n$. *Then* $\bar{a} \in K(\bar{G})$ *and* $\bar{G} = G/H$ *is a p-group, where H is the set of $p'$-elements.*

PROOF. This follows from Lemma 3.1.(iv) and Lemma 3.6 and its proof.

LEMMA 3.8. *If* $[h,a] = 1$, $h \in \eta(a)$ *and* $|a| = p^n$, *then* $ah \notin \eta(a)$ *if and only if* $(p, |<h>/(<h> \cap <a>)|) = 1$.

PROOF. Note first that $h \in \eta(a)$ implies that $<h> \cap <a> \leq <a^p>$. Let $m = |(<h> \cap <a>)|$. Then $h^m = (a^p)^s$ for some $s$. If $1 = p^n k + ml$, then $(ah)^{ml} = a^{ml} h^{ml} = a^{1-p^{sk}} a^{psl} = a^{psl+1}$. Since $(p, psl+1) = 1$, there exists a $t$ such that $(ah)^{mlt} = (a^{psl+1})^t = a$, so $ah \notin \eta(a)$. In the other direction, if $(ah)^q = a$, then $a^{q-1} = h^{-q} \in <h> \cap <a> \leq <a^p>$. This gives both $p \mid q-1$ and $m \mid q$, from which we conclude $(p,m) = 1$.

LEMMA 3.9. *Suppose* $g \in \eta(a)$ *and* $x^p = a$. *Then if* $g \notin \eta(y)$, $y \in \eta(x)$.

PROOF. By Lemma 3.2 (ii), $y \in \eta(a)$. Thus $x^p = a \notin <y>$, so $x \notin <y>$.

LEMMA 3.10. **[8; p.191]**. *If H is an equipotent subgroup of G, then* $K(G) \leq K(H)$.

DEFINITION. *Suppose G has a normal $p^\infty$-subgroup C of finite index. We say*

(i) *G is of type* $T_1$ *for the prime p if C is central;*

(ii) *G is of type* $T_2$ *if $p = 2$, C is not central, and every element x not in* $C_G(C)$ *satisfies* $x^{-1}cx = c^{-1}$ *for all* $c \in C$ *and there exists an m such that* $x^{2m} = a$, *the unique element in C of order 2.*

LEMMA 3.11. *The group G has an infinite abelian subgroup of finite index, and* $K(G) \neq E$ *if and only if one of the following holds:*

(i) *G is of type* $T_1$, *in this case* $K(G) = C$;

(ii) *G is of type* $T_2$, *in this case* $K(G) = <a>$.

PROOF. If G has an infinite abelian subgroup of finite index, then G has an infinite normal abelian subgroup $A$ of finite index. Now $K(G) \leq K(A) = E$ unless $A$ has a characteristic $p^\infty$-subgroup $C$ of finite index. Thus $C \triangleleft G$, and $G$ is a finite extension of $C$. Suppose $C$ is not central in $G$. Then $K = K(G) \leq C$ is central and finite and contains all the elements of $C$ which have order $p$. Let $x$ be any element in $G$ which fails to centralize $C$. Then conjugation by $x$, $\alpha_x$, is a non-trivial automorphism of $C$, which fixes every element of order $p$. By a theorem of Baer (see **[6; Lemma 3.28]**), $\alpha_x$ has infinite order unless $p = 2$. Thus we have a contradiction unless $p = 2$. For $p = 2$, we have the same contradiction if $\alpha_x$ fixes every element of order 4, so $\alpha_x^2 = 1$, $x^2 \in C_G(C)$ and $\alpha_x(c) = c^{-1}$ for every $c \in C$. Now $(xc)^2 = x^2(x^{-1}cx)c = x^2$ for each $c \in C$. Since $\eta(a)$ is finite where $a$ is the element of order 2 in $C$, $a \in <xc>$ for all but finitely many $c$. Suppose that there are infinitely many $c \in C$ such that $a = (xc)^{m_c}$ for $m_c$ odd. Then $(xc)^{m_c} = x^{m_c}c$, and since $x$ has finite order, there are infinitely many $m_c$ which are equal to $m$, some fixed integer. This gives infinitely many $c \in C$ such that $x^m c = a$ for a fixed integer $m$, clearly a contradiction. Thus for all but finitely many $c$, $a = (xc)^{m_c}$, with $m_c$ even. For these $c$, $a = (xc)^{m_c} = x^{m_c}$, so, in fact, there exists an even integer $2m$ such that $a = x^{2m}$. This shows that $G$ is type $T_2$. Clearly $K(G) = <a>$ since no other non-trivial element of C is central.

Conversely, if $G$ is of type $T_1$, $K(G) = C$ by Theorem C(iii). If $G$ is of type $T_2$, let us compute $\eta(a)$. Since each element $x$ not in $C_G(C)$ satisfies $x^{2m} = a$ for some $m$, these $x$ are not in $\eta(a)$. On the other hand, $C_G(C)$ is a finite central extension of a $2^\infty$-group, so $\eta(a)$ contains only finitely many elements in $C_G(C)$. This shows that $K(G) \geq <a>$. In addition, $K(G) \leq C$, and $a$ is the only non-trivial central element in $C$. Thus $K(G) = <a>$.

**§4. Main Results**. First we characterize the groups of type $T_1$ and type $T_2$ (see Lemma 3.11).

THEOREM 4.1.

(i) $G$ is of type $T_1$ for the prime $p$ if and only if $G = (H \times C)_A$, where $H$ is a finite group, $C$ is a $p^\infty$-group, and $A \leq Z(H)$ is a cyclic group of order $p^n$;

(ii) $G$ is of type $T_2$ if and only if $G = <x, G_1>$, where $G_1 = (H \times C)_A$ is of type $T_1$ for the prime 2, $x^{-1}cx = c^{-1}$ for every element $c$ of $C$, $[G : G_1] = 2$ and there exists an $m$ such that $x^{2m} = a$, the unique element of $C$ of order 2.

PROOF.

(i) If $G$ is of type $T_1$ then $G = \bigcup_{i \leq n} x_i C$. Let $H = <\{x_i : i \leq n\}>$. Since $G$ is locally finite, $H$ is finite. Clearly $G = HC$, $H \cap C \leq Z(G)$, and is cyclic.

(ii) If $G$ is of type $T_2$, then $G_1 = C_G(C)$ is of type $T_1$ for the prime 2. Suppose $x$ and $y \notin C_G(C)$. Then $x^{-1}cx = c^{-1} = y^{-1}cy$ for every $c \in C$, so $yx^{-1} \in C_G(C)$. Thus $[G : C_G(C)] = 2$ and $G = <x, C_G(C)>$ for any $x \notin C_G(C)$. In the other direction, note that for every $g \in G$ either $g^{-1}cg = c^{-1}$ for every $c \in C$ or $g \in C_G(C)$, depending upon the parity of the number of times that $x$ appears in the word which generates $g$. Again, we deduce $[G : C_G(C)] = 2$, so $G_1 = C_G(C)$.

EXAMPLE 4.2.

(i) Let $H$ be a non-abelian group of order $p^3$ and exponent $p$. Let $A$ be the center of $H$. Then $G = (H \times C)_A$ is not isomorphic to $C \times F$ for any finite group $F$ since $|G/C| = p^2$ and if, $|F| = p^2$, then $C \times F$ would be abelian.

(ii) Let $G_1 = (H \times C)_A$ be a group of type $T_1$. Let $y$ be any element in $G_1 - C$ such that $y^m = a$. Suppose there exists an automorphism $\alpha$ of $G_1$ which fixes $y$, $\alpha(c) = c^{-1}$ for each $c \in C$ and $\alpha^2(g_1) = y^{-1}g_1 y$ for every $g_1 \in G$. Then there exists a group $G$ such that $G = <x, G_1>$ with $G/G_1 \cong Z_2$, $x^2 = y$, $x^{-1}g_1 x = \alpha(g_1)$ for every $g_1 \in G$ (see [7, 9.7.1]). If $G_1$ is abelian and $m = 1$, then such an $\alpha$ exists. One can easily construct examples with $m \neq 1$. For example, let $G_1 = <z> \times C$, $z^2 = 1$. Then if $y = (z, a_2)$ with $a_2$ an element of order 4, we can define $\alpha$ by $\alpha(zc) = za_2^2 c^{-1}$ for each $c \in C$.

Next we consider countable groups with non-trivial finite $K$.

THEOREM 4.3. *If $G$ is a countable group with $K(G) \cong Z_{p^n}$, then $G$ has an infinite section $H$ of finite index with $K(H) = <a> \cong Z_p$, where $H$ is a $p$-group and $\eta(a, H) = E$.*

PROOF. First let us show that every group $G$ which is of type $T_2$ has a section $S$ of finite index which is infinite generalized quaternion, that is, $S = <x, C>$ with $C \cong Z_{2^\infty}$, $x^{-1}cx = c^{-1}$ for every element $c$ of $C$, and $x^2 = a$, the unique element of $C$ of order 2. Note that such a group S has $K(S) = <a>$ and $\eta(a) = E$. We shall induct on the integer $m$ in the definition of type $T_2$. Let $G = <x, G_1>$ be of type $T_2$ as in Theorem 4.1 (ii) with $G_1 = (H \times C)_A$, etc. Let $a_2$ be an element of $C$ of order 4. Let $z = x^m a_2$. We have $z^{-1}cz = x^{-m}cx^m = c^{(-1)^m}$ for every $c \in C$ and $z^2 = (x^m a_2)^2 = x^m a_2 x^{-m} x^{2m} a_2 = a_2^{(-1)^m} a_2 a$.

If $m$ is odd, $<z,C>$ is the desired section. If $m$ is even, say $m = 2k$, then $S_1 = <x,C>/N$ with $N = \{1, z, za, a\}$ has the desired generalized quaternion section by induction. To see this, note that we have shown that $(za)^2 = z^2 = 1$; $[z.c] = [x^m, c] = 1$ for every $c \in C$, $x^{-1}zx = xx^m a_2 x = x^m a_2^{-1} = za$ so that $N \triangleleft S_1$; $a_2 N$ is the unique element in $CN/N$ of order 2 since $a_2 \notin N$; and finally, $x^{2k} = x^m = za_2^{-1} = zaa_2$, so $x^{2k}N = a_2 N$.

Let $K(G) = <a> \cong Z_{p^n}$. We may assume that $G$ does not have a subgroup $H_1$ of finite index which is type $T_2$, otherwise $H_1$ has a subgroup $H$ which satisfies the conditions of the theorem by the argument just given.

Since by Lemma 3.3, $\alpha(\eta(a)) = \eta(a)$ for every $\alpha \in Aut(G)$, and since $\eta(a)$ is finite, each element in $\eta(a)$ has finitely many conjugates. Thus $H = C(\eta(a))$ has finite index in $G$. If $C = K(H)$ were a $p^\infty$-group, then since $C$ has finite index in $H$ by Theorem C(iii), $C$ would have finite index in $G$. Since $K(G)$ is finite, $C$ could not be central. Thus by Lemma 3.11, we would have $K(G) = E$; a contradiction. Thus $K(H)$ is finite. Since by Lemma 3.10, $K(G) \leq K(H)$, $K(H)$ is non-trivial, and we may suppose $G = C(\eta(a))$. The set $N$ of $p'$-elements forms an abelian subgroup and $N \leq \eta(a)$ by Lemma 3.6. Now $\bar{a} \in K(G/N)$ by Lemma 3.7. If $P/N = K(G/N) \cong Z_{p^\infty}$, then $[G/N : P/N] < \infty$. Now by **[7,3.2.10]**, $P$ is abelian, hence $P = Q \times N$, where Q is a $p^\infty$-group and $N$ is a finite $p'$-group. Thus $[G : Q] < \infty$, so $K(G) \cong Z_{p^\infty}$; a contradiction. Thus we may assume that G is a $p$-group.

We have shown in Lemma 3.5 that $K(G/<a^p>) = <\bar{a}>$, so we may assume that $K(G) \cong Z_p$. Consider $Z(G)$. If it is infinite, then by Theorem A, $K(Z(G)) = E$ unless $Z(G) = Q \times F$, where $Q$ is a $p^\infty$-group and $F$ is finite. If $Z(G)$ is finite, $Z(G) = Q \times F$, where $Q$ is cyclic and $F$ is finite. Thus $Z(G) = Q \times F$, where $Q$ is cyclic or quasicyclic and $a \in Q$. Now $\eta(a) = \{qf : q \in Q, f \in F, |q| \leq |f|\}$.

.

Let us compute in $\bar{G} = G/F$. Suppose $|\eta(\bar{x})| < \infty$, with $x \notin Z(G)$. Then $|G - C(x)|$ is infinite and so is $|(G - C(x))/F|$. Since $C(xf) = C(x)$ for all $f \in F$, $G - C(x) = G - C(xf) \leq \eta(xf)$. If $y \in G - C(x)$, but $y \notin \eta(\bar{x})$, then $\bar{x} = \bar{y}^n$ or $y^n = xf$ for some $n$ and $f \in F$, contradicting $y \in \eta(xf)$. This shows that $|\eta(\bar{x})|$ is finite only if $x \in Z(G)$.

.

Let us show that $\eta(\bar{a}) = \{\bar{1}\}$. Suppose $\bar{z}^p = \bar{1}$, $\bar{z} \neq \bar{1}$, If $z \in \eta(a)$, then $z = qf$, with $q \in Q, f \in F, |q| \leq |f|$. Since $q^p f^p = z^p \in F$, $q^p = 1$, so $q \in <a>$. Thus $<\bar{z}> = <\bar{a}>$. If $z \notin \eta(a)$, then since $z^p \in F \leq \eta(a)$, $z^p = 1$. Since $z \notin F$, $<z> = <a>$ and again $<\bar{z}> = <\bar{a}>$. Thus $<\bar{a}>$ is the only subgroup of $G$ of order $p$. It follows that $\eta(\bar{a}) = \{\bar{1}\}$.

If $\eta(\bar{x})$ is finite, then $x \in Z(G)$, so $x = qf$ where $q \in Q$. We shall show that $\eta(\bar{q}) = <\bar{q}^p>$. We have shown that $\eta(\bar{a}) = \{\bar{1}\}$. Suppose that we have established the equality $\eta(\bar{q}) = <\bar{q}^p>$ for all $q \in Q$ with $[<q>:<a>] < p^m$, and there is a $q \in Q$, with $q^{p^m} = a$. Suppose $\bar{y} \in \eta(\bar{q}) - <\bar{q}^p>$. Then since $[<q^p>:<a>] < p^m$, $\bar{y} \notin \eta(q^p)$, so $\bar{y}^n = \bar{q}^p$ for some $n$. Since $\bar{y} \notin <\bar{q}^p>$, we may assume that $p | n$, say, $n = pl$. Thus $\bar{z}^p = (\bar{y}^l)^p = \bar{q}^p$, where $\bar{z} = \bar{y}^l \in \eta(\bar{q}) - <\bar{q}^p>$. We have $(\overline{zq^{-1}})^p = 1$. Since $\bar{y} \in \eta(\bar{q})$, $\bar{z} \neq \bar{q}$. Since $\eta(\bar{a}) = \{\bar{1}\}$, $\overline{zq^{-1}} = \bar{a}^s$ or $\bar{z} = \overline{qa^s} \in <\bar{q}>$. This contradicts $\bar{z} \in \eta(\bar{q}) - <\bar{q}^p>$. We have shown that $K(G) \cong Q$ and that $\eta(\bar{q}) = <\bar{q}^p>$ for every non-trivial $q \in Q$. If $Q$ is a $p^\infty$-group, then $[G:K(G)] < \infty$, so $[G:Q] < \infty$, which gives $K(G) \cong Z_{p^\infty}$, which is a contradiction. Thus there is a group G with $K(G) = <a>$ and $\eta(x) = <x^p>$ for each non-trivial $x \in <a>$. Now in $\bar{G} = G/<a^p>$, $K(\bar{G}) = <\bar{a}> \cong Z_p$ by Lemma 3.5 and $\eta(\bar{a}) \leq \bar{\eta}(a) = <\bar{a}^p> = \{\bar{1}\}$. This proves the theorem.

We need the following lemma which is stated in **[10]**. The main ideas for the proof can also be found in **[6; p. 70, part 1]** and **[2; l.G.4 and l.G.6]**.

LEMMA 4.4. *Let G be infinite and either binary finite or a 2-group. If G has a finite maximal elementary abelian subgroup, then it has an infinite abelian subgroup of finite index.*

THEOREM 4.5. *If G is a countable group such that $K(G) \neq E$, then one of the following holds:*

*(i) G is of type $T_1$: $G = (P \times F)_C$, with F a finite group; P a $p^\infty$-group; and $C \leq Z(F)$, a cyclic group of order a power of p. In this case $K(G) = P$.*

*(ii) G is of type $T_2$: $G = <x, G_1>$, where $G_1 = (P \times F)_C$ with F a finite group; P a $2^\infty$-group; $C \leq Z(F)$, a cyclic group of order a power of 2, $x^{-1}zx = z^{-1}$ for every element z of P, $[G:G_1] = 2$ and there exists an m such that $x^{2m} = a$, the unique element of C of order 2. In this case $K(G) = <a>$.*

*(iii) G has an infinite section H, which is a 2-generated p-group for an odd prime p and $K(H) \cong Z_p$.*

PROOF. We know by Theorems C(iii) and 4.1 that $K(G) \cong Z_{p^\infty}$ if and only if (i) holds. On the other hand, if $K(G) \cong Z_{p^n}$ then by Theorem 4.3 there is an infinite section $H$ with $K(H) = <a> \cong Z_p$, with $H$ a $p$-group and $\eta(a) = \{1\}$. Clearly $<a>$ is a maximal elementary abelian subgroup. If $H$ is binary finite or a 2-group, then by Lemma 4.4, $H$ has an infinite abelian subgroup of finite index. By Lemma 3.11, this is only possible if (ii) holds. Otherwise $p \neq 2$ and $H$ has an infinite subgroup $L$ which is 2-generated with $a \in L$ and $K(L)$ finite. By Lemma 3.5, $L$ has a factor group which satisfies (iii).

COROLLARY 4.6. *If G is countable and either binary finite or a 2-group such that $K(G) \neq E$, then either G is of type $T_1$ and K(G) is $p^\infty$-group or G is of type $T_2$ and $K(G) \cong Z_2$.*

REMARKS. We have not been able to improve Theorem 4.5 to include in (iii) the statement $K(H) = <a>$ with $\eta(a) = \{1\}$. Of course, the likely place to look for groups of the type described in 4.5 (iii) is a central extension of $Z_p$ by a group $B$ which is either a Novikov-Adjan group (see **[4]**) or a Tarski-Monster constructed by Ol'shanskii **[5]**. This necessitates the study of maps $W : B \times B \to Z_p$ which satisfy

$$W(xy, z) + W(x, y) = W(x, yz) + W(y, z),$$

which seems quite difficult.

Now we examine some uncountable groups.

THEOREM 4.7. *If G is an uncountable class 2 nilpotent p-group for p odd, then $|\Omega_1(G)| = |G|$. For each infinite cardinal m there is a class 2 nilpotent 2-group G of power $2^m$ such that $|\Omega_1(G)| = 2^{<m}$.*

PROOF. Consider the tree $T = 2^{<m}$ of functions from ordinals less than $m$ into 2, ordered by function extension. $T$ has $2^{<m}$ nodes and $2^m$ paths (a path corresponds to a function from $m$ into 2.) Let $V$ be the $2^m$-dimensional vector space over $Z_2$ with basis $\Sigma$, the set of paths of $T$. Let $W$ be the $2^{<m}$-dimensional vector space over $Z_2$ with basis $N$, the set of nodes of $T$. We define for each $f$ and $g$ in $\Sigma$, $\rho(f, g)$ to be the largest element of $N$ common to both $f$ and $g$ if $f \neq g$ and $\rho(f, f) = 0$. Extend $\rho$ by bilinearity to all of $V$. Let $\gamma$ be any bilinear function from $V \times V$ to $W$ such that $\gamma(v, v) \neq 0$ for $v \in \Sigma$ and $\rho(u, v) = \gamma(u, v) + \gamma(v, u)$. Consider the class 2 nilpotent 2-group $G = V\gamma W$ (see **[1]**). If $(v, a)^2 = 1$ with $0 \neq v \in V$ and $a \in W$, then $\gamma(v, v) = 0$. Suppose $v = \sum_{i=1}^{n} f_i$ with $f_i \neq f_j \in \Sigma$. Then $\gamma(v,v) = \sum_{i,j=1}^{n} \gamma(f_i, f_j) = \sum_{i<j} \rho(f_i, f_j)$. If $a$ is maximal in $T$ among the $\rho(f_i, f_j)$ then $a$ is unpaired among the values $\rho(f_i, f_j)$, so $\gamma(v,v) \neq 0$. Thus only the elements of $W$ have order 2. This proves the second statement.

If G is an uncountable class 2 nilpotent $p$-group for $p$ odd with $|\Omega_1(G)| < |G|$, then since the derived group G' is abelian, $|G'| < |G|$ (see **[8; p. 184, Corollary 1]**). Assume that $G$ is an example with smallest cardinality. It is clear that $|G|$ is a successor cardinal. Since $G/G'$ is abelian, it has an elementary abelian subgroup $H/G'$ with $|H/G'| = |G/G'| = |G|$. If $x, y \in H$, then $[x,y] \in Z(G)$ and $x^p \in G' \leq Z(G)$, so $[x,y]^p = [x^p, y] = 1$, and hence $\exp H' = p$. Let $\{x^\alpha H'\}$ be a basis for $H/H'$. Since $x_\alpha^p \in H'$, there is a set $T$ of power $|G|$ and an $a \in H'$ such that $x_\alpha^p = a$ for every $\alpha \in T$. But then $(x_\beta^{-1} x_\alpha)^p = [x_\alpha, x_\beta^{-1}]^{\binom{p}{2}} = 1$, so $|\Omega_1(G)| = |G|$.

COROLLARY 4.8. *For each infinite cardinal m there is a class 2 nilpotent 2-group G of power $2^m$ such that every abelian subgroup A satisfies $|A| \leq 2^{<m}$.*

PROOF. This is clear since if $|A| > 2^{<m}$, then $|\Omega_1(G)| \geq |\Omega_1(A)| = |A| > 2^{<m}$; a contradiction.

REMARK. The uncountable group constructed in **[3]** is similar.

COROLLARY 4.9. *If G is an uncountable class 2 nilpotent group, then $K(G) = E$.*

PROOF. If G is a counter-example of smallest cardinality and $1 \neq a \in K(G)$, then $|\eta(a)| < |G|$, so we may assume $|\eta(a)| \leq m$ and $|G| = m^+$. Thus, by arguments like those used in the proof of the theorem, we need only prove the corollary for groups $G$ such that $G/G'$ and $G'$ have exponent $p$. If $p$ is odd, then Theorem 4.7 gives a contradiction, so we suppose that $p = 2$. We have $a^2 = 1$, and if $x \in G - \eta(a) - \{a\}$, then $x^2 = a$. Let $\{x_\alpha\}$ be such that $\{x_\alpha G'\}$ is a basis for $G/G'$. Let $x$ be any fixed element of $T = \{x_\alpha\} - \eta(a)$. If $x_\alpha \in T$, then $(xx_\alpha)^2 = x^2 x_\alpha^2 [x, x_\alpha] = [x, x_\alpha]$. Thus there must be $S_0 \leq T$ such that $|S_0| = |T|$ and $[x, x_\alpha] = a$ for all $x_\alpha \in S_0$. Fix $y \neq x \in S_0$. Again there is $S_1 \leq S_0$ such that $|S_1| = |S_0|$ and $[y, x_\alpha] = a$ for all $x_\alpha \in S_1$. Fix $z \neq x, y \in S_1$ and find $S_2 \leq S_1$ such that $|S_2| = |S_1|$ and $[z, x_\alpha] = a$ for all $x_\alpha \in S_2$. Thus for each $x_\alpha \in S_2$, $(xyzx_\alpha)^2 = (xy)^2(zx_\alpha)^2[xy, zx_\alpha] = [x,y][z, x_\alpha][x,z][y,z][x, x_\alpha][y, x_\alpha] = a^6 = 1$. But $xyzx_\alpha$ are distinct members of $\eta(a)$ for $x_\alpha \in S_2$, which is a contradiction.

REMARK. Although it seems quite likely that Corollary 4.9 holds for all uncountable locally nilpotent groups, we have not even been able to extend it to class 3 nilpotent groups.

EXAMPLE 4.10. There are class 2 nilpotent $p$-groups with $D > K = E$. To see this, consider one of the groups G constructed (using G.C.H.) in **[1]**. G has the properties ($p$ is a prime):

(i) $|G| = \lambda^+$, where $\lambda$ is an infinite cardinal;

(ii) $G' = Z(G) \cong Z_p$;

(iii) if $A$ is an abelian subgroup of $G$, $|A| \leq \lambda$;

(iv) $\exp(G/G') = p$.

Thus if $G_\alpha \leq G$, where $|G_\alpha| < |G|$, $G'_\alpha = G'$. It follows that $D = G'$. On the other hand, we know that $K \leq G' = Z(G)$. Let $1 \neq a \in G'$. Since $\eta(a) = \{g \in G : g^p = 1\} - (<a> - \{1\})$, if $|\eta(a)| \leq \lambda$, take $x_0 \in G - \eta(a)$, $x_0^p = a$. Then $|C(x_0) - \eta(a) = \lambda^+|$, and for every $x \in C(x_0) - \eta(a)$, $x^p = a$. Then $(x^{-1}x_0)^p = x^{-p}x_0^p = 1$, so $x^{-1}x_0 \in \eta(a)$. Thus $|\eta(a)| = \lambda^+$, and so $K = E$.

## ACKNOWLEDGMENT

I wish to express my appreciation to Professor Scott for pointing out several errors in a previous version of this paper.## REFERENCES

1. A. Ehrenfeucht and V. Faber, *Do infinite nilpotent groups always have equipotent Abelian subgroups?* Kon. Nederl. Akad. Wet. A75 (1972), 202-209.

2. O. H. Kegel and B. A. F. Wehrfritz, *Locally Finite Groups*. North Holland (1973).

3. L. G. Kovacs, B. H. Neumann, and H. de Vries, Proc. Royal Soc. London A260 (1961),304-316.

4. P. S. Novikov and S. 1. Adjan, *On infinite periodic groups* (Russian), Izv. Akad. Nauk SSSR Ser. Mat. 32 (1968), 212-244, 251-524, 709-731. MR 29#1532 a, b, c.

5. A. Yu. Ol'shanskii, *Groups of bounded period with subgroups of prime order*, Algebra and Logic 21 (1982), 369-418.

6. D. J. S. Robinson, *Finiteness Conditions and Generalized Soluble Groups*, Parts 1 and 2, Ergebnisse der Math. Bde. 62/63. Springer-Verlag (1972).

7. W. R. Scott, *Group Theory*, Prentice Hall, Englewood Cliffs (1964).12